\documentclass{amsart}
\usepackage{amsmath, amsthm, amssymb, amscd, mathrsfs, eucal, palatino}

\begin{document}

\newtheorem{theorem}{Theorem}[section]
\newtheorem{lemma}[theorem]{Lemma}
\newtheorem{corollary}[theorem]{Corollary}
\newtheorem{conjecture}[theorem]{Conjecture}
\newtheorem{question}[theorem]{Question}
\newtheorem{problem}[theorem]{Problem}
\newtheorem*{claim}{Claim}
\newtheorem*{criterion}{Criterion}
\newtheorem*{main_thm}{Theorem A}

\theoremstyle{definition}
\newtheorem{definition}[theorem]{Definition}
\newtheorem{construction}[theorem]{Construction}
\newtheorem{notation}[theorem]{Notation}

\theoremstyle{remark}
\newtheorem{remark}[theorem]{Remark}
\newtheorem{example}[theorem]{Example}

\def\area{\text{area}}
\def\id{\text{id}}
\def\H{\mathbb H}
\def\Z{\mathbb Z}
\def\R{\mathbb R}
\def\Q{\mathbb Q}
\def\F{\mathcal F}
\def\SL{\textnormal{SL}}
\def\PSL{\textnormal{PSL}}
\def\fix{\textnormal{fix}}
\def\RP{\mathbb RP}
\def\PL{\textnormal{PL}}
\def\homeo{\textnormal{Homeo}}
\def\inte{\textnormal{int}}

\def\CAT{\textnormal{CAT}}
\def\Aut{\textnormal{Aut}}
\def\Out{\textnormal{Out}}
\def\Sol{\textnormal{Sol}}
\def\Nil{\textnormal{Nil}}
\def\til{\widetilde}

\title{Stable commutator length in subgroups of $\PL^+(I)$}
\author{Danny Calegari}

\begin{abstract}
Let $G$ be a subgroup of $\PL^+(I)$. Then the stable
commutator length of every element of $[G,G]$ is zero.
\end{abstract}

\address{Department of Mathematics \\ California Institute of Technology \\
Pasadena CA, 91125}
\email{dannyc@its.caltech.edu}
\date{8/13/2006, Version 0.04}

\maketitle

\section{Introduction}

The purpose of this note is to prove a vanishing theorem for stable commutator length
in groups of $\PL$ homeomorphisms of the interval. For convenience, we restrict attention
to subgroups of the group of orientation preserving $\PL$ homeomorphisms, denoted
in the sequel by $\PL^+(I)$, where $I$ denotes the unit interval $[0,1]$.
By a theorem of Bavard (see \S~\ref{background_section}), vanishing of stable
commutator length is equivalent to the injectivity of the map from bounded 
to ordinary cohomology in dimension $2$.  

Since the dimension of the second bounded cohomology of a
nonabelian free group is uncountable, 
this gives a new proof of the celebrated result of Brin--Squier (\cite{Brin_Squier})
that $\PL^+(I)$ does not contain a nonabelian free subgroup.

There are at least two other important classes of groups for which vanishing of
stable commutator length is known to hold:

\begin{itemize}
\item{Irreducible lattices in semisimple Lie groups of rank at least
$2$; this follows from more general work of Burger and Monod \cite{Burger_Monod}.}
\item{Amenable groups; in this case, bounded cohomology with real coefficients
vanishes in every dimension by Trauber's theorem (see \S~\ref{background_section}) and
therefore the map is injective for trivial reasons.}
\end{itemize}

It is considered an important open question whether or not Thompson's group
$F < \PL^+(I)$, which consists of homeomorphisms with dyadic rational slopes and break points,
is amenable. More generally, no counterexamples are known to the conjecture that
a finitely presented torsion free group with the property that
every subgroup has vanishing stable commutator length is amenable
(this should perhaps be thought of as a kind of ``homological" version of
von Neumann's conjecture). These and related problems are one of the main motivations for
this paper.

\subsection{Acknowledgements}
While writing this paper, I was partially supported by a Sloan Research Fellowship, and
NSF grant DMS-0405491. I also benefited from useful discussions with \'Etienne Ghys and
Dieter Kotschick. I would further like to thank the anonymous referee, many of whose comments have
been incorporated verbatim into this paper.

\section{Background material}\label{background_section}

\begin{definition}\label{bounded_cohomology}
Let $G$ be a group, and $C_*(G)$ the (bar) complex of integral $G$-chains.
Let $C^*(G)\otimes \R$ be the dual complex of real-valued cochains.
For each $n$, let $C^n_b(G)\otimes \R$ denote the vector space of cochains $f$
for which $\sup_\sigma |f(\sigma)|$ is finite, where $\sigma$ ranges over the
generators of $C_*(G)$. The (real) {\em bounded cohomology} of $G$, denoted
$H^*_b(G;\R)$, is the cohomology of the complex $C^*_b(G)\otimes \R$.
\end{definition}

Note that $H^n_b(G;\R)$ carries an $L^\infty$ pseudo-norm for each $n$.

Bounded cohomology behaves well under amenable extensions:

\begin{theorem}[Trauber]\label{Trauber_theorem}
Let $$1 \to H \to G \to A \to 1$$
be a short exact sequence of groups, where $A$ is amenable. Then
the natural homomorphisms $H^*_b(G;\R) \to H^*_b(H;\R)$ are isometric
injections.
\end{theorem}

See e.g. \cite{Gromov_volume} page 39 for a proof.

\begin{definition}
Let $G$ be a group, and $[G,G]$ the commutator subgroup. For any
$g \in [G,G]$, the {\em commutator length} of $g$, denoted
$\ell'(g)$, is the minimal number of commutators whose product is equal
to $g$. The {\em stable commutator length}, denoted $\ell(g)$ is defined
to be
$$\ell(g) = \liminf_{n \to \infty} \frac {\ell'(g^n)} {n}$$
\end{definition}

By including bounded cochains in all cochains, one obtains 
canonical homomorphisms
from bounded cohomology to ordinary cohomology.
There is a fundamental relationship between stable commutator
length and bounded cohomology, discovered by Bavard:

\begin{theorem}[Bavard]\label{Bavard_theorem_1}
Let $G$ be a group. Then the canonical map from
bounded cohomology to ordinary cohomology $H^2_b(G;\R) \to H^2(G;\R)$
is injective if and only if the stable commutator length vanishes on $[G,G]$.
\end{theorem}

See \cite{Bavard}.

Bavard's theorem makes use of the notion of quasimorphisms:

\begin{definition}
Let $G$ be a group. A (homogeneous) {\em quasimorphism} on $G$ is a map
$$f:G \to \R$$
for which there is some smallest $\epsilon(f) \ge 0$ (called the {\em error}
or {\em defect} of $f$) such that
$$f(a^n) = nf(a)$$
and
$$|f(a) + f(b) - f(ab)| \le \epsilon(f)$$
for all $a,b \in G$.
\end{definition}
Note that a homogeneous quasimorphism is necessarily a class function.

The set of all homogeneous quasimorphisms on $G$, denoted $Q(G)$, 
has the structure of a vector
space. Quasimorphisms with error $0$ are homomorphisms. There is an
exact sequence

$$1 \to H^1(G;\R) \to Q(G) \xrightarrow{\delta} H^2_b(G;\R) \to H^2(G;\R)$$
where $\delta$ denotes the coboundary map. See \cite{Bavard} for a proof.

Thus Bavard's theorem may be interpreted as saying that if $G$ is a group,
the quotient $Q(G)/H^1(G;\R)$ is zero exactly when the
stable commutator length of every element of $[G,G]$ vanishes.

In terms of quasimorphisms, Bavard proves the following sharper statement:

\begin{theorem}[Bavard]\label{Bavard_theorem_2}
Let $G$ be a group, and $g \in [G,G]$. Then
$$\ell(g) = \frac 1 2 \sup_{f \in Q(G)/H^1(G)} \frac {\|f(g)\|} {\epsilon(f)}$$
\end{theorem}

Theorem~\ref{Bavard_theorem_1} and Theorem~\ref{Trauber_theorem} together
imply that $\ell(g)=0$ for any $g \in [G,G]$ whenever $G$ is an amenable group.

\section{Subgroups of $\PL^+(I)$}

Given a subgroup $G < \PL^+(I)$ we denote by $\fix(G)$ the set of 
common fixed points of all elements of $G$. 

\begin{definition}
The {\em endpoint homomorphism} is the homomorphism
$$\eta: \PL^+(I) \to \R \oplus \R$$
defined by
$$\eta(g) = (\log g'(0), \log g'(1))$$
Given $G < \PL^+(I)$, denote by $G_0$ the kernel of $\eta$ restricted
to $G$.
\end{definition}

Observe that every $g \in G_0$ fixes a neighborhood of both $0$ and $1$.

\begin{main_thm}
Let $G$ be a subgroup of $\PL^+(I)$. Then the stable
commutator length of every element of $[G,G]$ is zero.
\end{main_thm}
\begin{proof}
{\noindent \bf Case 1:} $\fix(G) = \lbrace 0,1 \rbrace$.

Let $G_0$ be the kernel of $\eta:G \to \R \oplus \R$.
Let $K = [G_0,G_0]$ and let $g \in K$. Then we can write
$$g = [a_1,b_1][a_2,b_2] \cdots [a_m,b_m]$$
for some integer $m$ and $a_i,b_i$ in $G_0$. Let $J$ be the smallest interval which
contains the support of all the $a_i,b_i$ and $g$. Then $J$ is properly contained in
$(0,1)$. Since $\fix(G)$ contains no interior points, 
there is some $j \in G$ with $j(J) \cap J = \emptyset$, and therefore $j^n(J) \cap J = \emptyset$
for all nonzero $n$.

Let $G_0(J)$ be the subgroup of $G_0$ consisting of elements with support
contained in $J$. For each $n$ we define a diagonal monomorphism
$$\Delta_n: G_0(J) \to G_0$$ by
$$\Delta_n(c) = \prod_{i=0}^n c^{j^i}$$
where the superscript notation denotes conjugation.
Define
$$g' = \prod_{i=0}^n (g^{i+1})^{j^i}$$
Then 
$$[g',j]=\Delta_n(g) (g^{-n-1})^{j^{n+1}}$$
On the other hand, 
$$\Delta_n(g) = \Delta_n([a_1,b_1][a_2,b_2] \cdots [a_m,b_m]) = 
[\Delta_n(a_1),\Delta_n(b_1)] \cdots [\Delta_n(a_m),\Delta_n(b_m)]$$
and therefore $g^{n+1}$ can be written as a product of at most $m+1$ 
commutators in elements of $G$. Since $m$ is fixed but $n$ is arbitrary, 
it follows that the stable commutator length
of $g$ is zero, and hence $f([G_0,G_0])=0$ for every quasimorphism $f \in Q(G)/H^1(G)$.

Now, let $g \in [G,G]$. Observe that $[G_0,G_0]$ is normal in $G$, so we
can form the quotient $H = G/[G_0,G_0]$ which is two-step solvable, and therefore
amenable. Let $\varphi:G \to H$ be the quotient homomorphism.
By Trauber's Theorem~\ref{Trauber_theorem} and Bavard's Theorem~\ref{Bavard_theorem_1}, 
$\ell(\varphi(g)) = 0$ in $H$. This means
that we can write
$$g^n = [a_1,b_1]\cdots [a_m,b_m] c$$
where $c \in [G_0,G_0]$, where $n$ is arbitrarily  big, and 
$m/n$ is as small as we like. Let $f$ be a quasimorphism of defect at most $1$.
By the above, we have $f(c) = 0$, and therefore $f(g^n) \le 2m+1$, and
$f(g) \le \frac {2m+1} {n}$. Since $n$ is arbitrarily big, and $m/n$ is
as small as we like, $f(g)=0$. Since $f$ and $g$ were arbitrary, $Q(G)/H^1(G)=0$.

Applying Theorem~\ref{Bavard_theorem_1},
this proves the theorem when $\fix(G) = \lbrace 0,1 \rbrace$.

\vskip 12pt

{\noindent \bf Case 2:} $\fix(G)$ is arbitrary.

Suppose $f \in Q(G)$ has defect at most $1$, and suppose
$f(g) \ne 0$ where $g \in [G,G]$. Let $H$ be a finitely generated subgroup
of $G$ such that $g \in [H,H]$.
Then $\fix(H)$ is equal to the intersection of $\fix(h_i)$ for the generators $h_i$.
The fixed set of any element of $\PL^+(I)$ is a union of finitely many
points and intervals, so the same is true for $\fix(H)$.
Hence $I \backslash \fix(H)$ consists of finitely many open intervals, whose
closures we denote by $I_1,I_2, \dots I_n$. 

Let $\rho: H \to \R^{2n}$ denote the product of the endpoint homomorphisms for
each $i$, and let $H_0$ denote the kernel. We will show that $f$ vanishes on $[H,H]$,
contrary to the fact that $f(g) \ne 0$ and $g \in [H,H]$.

Let $r \in [H_0,H_0]$ and suppose we have an expression
$$r = [a_1,b_1] \cdots [a_m,b_m]$$
where all the $r,a_i,b_i$ have support in the union $\cup_i J_i$ where
$J_i \subset \inte(I_i)$ is an interval for each $i$. For each $i$ there is $j_i \in H$
with $j_i(J_i) \cap J_i = \emptyset$. Note that this implies $j_i^n(J_i) \cap J_i = \emptyset$
for all nonzero $n$.

However, we claim that we can construct 
a {\em single element} $j \in H$ such that $j(J_i) \cap J_i = \emptyset$ 
for all $i$ {\em simultaneously}.

\vskip 12pt

The case $n=1$ is trivial; in the interests of exposition we describe the situation 
$n=2$ in detail before moving onto the general case.

Without loss of generality, we may assume $j_1$ moves $J_1$ to the right.
Now, let $J_2'$ be the smallest interval which contains both $J_2$ and $j_1^{-1}(J_2)$,
and let $j_2$ be such that $j_2(J_2)' \cap J_2' = \emptyset$.
After replacing $j_2$ by $j_2^{-1}$ if necessary, we may also assume that $j_2$
moves the leftmost point of $J_1$, which we denote by $J_1^-$, to the right. We
also use the notation $J_1^+$ to denote the rightmost point of $J_1$.
i.e. 
$$j_1(J_1^-) > J_1^+, \; j_2(J_1^-) \ge J_1^-$$
Then 
$$j_1j_2(J_1^-) > J_1^+$$
and therefore 
$$j_1j_2(J_1) \cap J_1 = \emptyset$$
Moreover,
$$j_1j_2(J_2) \cap J_2 = j_1(j_2(J_2) \cap j_1^{-1}(J_2)) \subset j_1(j_2(J_2') \cap J_2') = \emptyset$$

Now we treat the general case. As before, without loss of generality, we
assume $j_1$ moves $J_1$ to the right. For all $i > 1$ we let $J_i' \subset I_i$ denote
the smallest interval which contains both $J_i$ and $j_1^{-1}(J_i)$. By
induction, we assume that there is some $j$ with $j(J_i') \cap J_i' = \emptyset$
for all $i > 1$ simultaneously. After replacing $j$ with $j^{-1}$
if necessary, we may assume that $j$ moves the leftmost point of $J_1$ to the right.
Then the argument above shows that 
$$j_1j(J_i) \cap J_i = \emptyset$$
for all $i$. Therefore we have proved the claim.

\vskip 12pt

But now the proof that $\ell(r)=0$ follows exactly as in Case 1, since for any $m$
there is a diagonal monomorphism 
$$\Delta_m:G_0(\cup_i J_i) \to G_0$$
defined by
$$\Delta_m(c) = \prod_{i=0}^m c^{j^i}$$
where now $j$ moves every $J_i$ off itself simultaneously. Since $r$ was
arbitrary, it follows that stable commutator length vanishes on $[H_0,H_0]$, 
and since $H/[H_0,H_0]$ is amenable, $f$
must vanish on all of $[H,H]$ by Theorem~\ref{Bavard_theorem_1}, 
contrary to the definition of $H$.
This contradiction implies
that $Q(G)/H^1(G)=0$, and the theorem follows.
\end{proof}

\begin{remark}
The ``diagonal trick" is a variation on Mather's argument (\cite{Mather}) to
prove the acyclicity of $\homeo_0(\R^n)$. This argument was modified
by Matsumoto-Morita (\cite{Matsumoto_Morita}) to prove vanishing of
all the bounded cohomology of $\homeo_0(\R^n)$. One significant difference
between $\PL^+(I)$ and $\homeo_0(\R^n)$ is that
every finitely generated subgroup $G$ of $\homeo_0(\R^n)$ is contained
in an {\em unrestricted} wreath product with $\Z$ (i.e. a product
of the form $\prod_i G \rtimes \Z$), whereas in a $\PL$
group, only {\em restricted} wreath products (i.e. products of
the form $\oplus_i G \rtimes \Z)$ with infinite groups are possible.
\end{remark}

\vfill
\pagebreak

\end{document}